\magnification 1250
\pretolerance=500 \tolerance=1000  \brokenpenalty=5000
\mathcode`A="7041 \mathcode`B="7042 \mathcode`C="7043
\mathcode`D="7044 \mathcode`E="7045 \mathcode`F="7046
\mathcode`G="7047 \mathcode`H="7048 \mathcode`I="7049
\mathcode`J="704A \mathcode`K="704B \mathcode`L="704C
\mathcode`M="704D \mathcode`N="704E \mathcode`O="704F
\mathcode`P="7050 \mathcode`Q="7051 \mathcode`R="7052
\mathcode`S="7053 \mathcode`T="7054 \mathcode`U="7055
\mathcode`V="7056 \mathcode`W="7057 \mathcode`X="7058
\mathcode`Y="7059 \mathcode`Z="705A
\def\spacedmath#1{\def\packedmath##1${\bgroup\mathsurround
=0pt##1\egroup$}
\mathsurround#1
\everymath={\packedmath}\everydisplay={\mathsurround=0pt}}
\def\nospacedmath{\mathsurround=0pt
\everymath={}\everydisplay={} } \spacedmath{2pt}

\def\phfl#1#2{\normalbaselines{\baselineskip=0pt
\lineskip=10truept\lineskiplimit=1truept}\nospacedmath\smash 
{\mathop{\hbox to 8truemm{\rightarrowfill}}
\limits^{\scriptstyle#1}_{\scriptstyle#2}}}
\def\hfl#1#2{\normalbaselines{\baselineskip=0truept
\lineskip=10truept\lineskiplimit=1truept}\nospacedmath\smash
{\mathop{\hbox to
10truemm{\rightarrowfill}}\limits^{\scriptstyle#1}_{\scriptstyle#2}}}
\def\diagram#1{\def\normalbaselines{\baselineskip=0truept
\lineskip=10truept\lineskiplimit=1truept}   \matrix{#1}}
\def\vfl#1#2{\llap{$\scriptstyle#1$}\left\downarrow\vbox to
5truemm{}\right.\rlap{$\scriptstyle#2$}}

\def\iso{\vbox{\hbox to .8cm{\hfill{$\scriptstyle\sim$}\hfill}
\nointerlineskip\hbox to .8cm{{\hfill$\longrightarrow $\hfill}} }}

\font\eightrm=cmr8     \font\sixrm=cmr6  
\def\pc#1{\tenrm#1\sevenrm}
\def\tx{\kern-1.5pt -}
\def\cqfd{\kern 2truemm\unskip\penalty 500\vrule height 4pt depth 0pt 
width 4pt\medbreak} 
\def\no{n\up{o}\kern 2pt}
\def\ind{\par\hskip 1truecm\relax}
\def\indp{\par\hskip 0.5truecm\relax}
\def\Pic{\mathop{\rm Pic}\nolimits}
\frenchspacing
\vsize = 25truecm
\hsize = 16truecm
\voffset = -.5truecm
\parindent=0cm
\baselineskip15pt

\centerline{\bf  On the Chow ring of a K3 surface}
\smallskip
\smallskip \centerline{Arnaud {\pc BEAUVILLE}}
\vskip1.2cm 

{\bf Introduction}
\smallskip
\ind   An important algebraic invariant
of a projective manifold $X$ is the  Chow ring $CH(X)$  of algebraic cycles
on $X$ modulo linear equivalence. It is  graded by the codimension  of
cycles; the ring structure comes from the  intersection product.
 For a surface we have 
$$CH(X)={\bf Z}\oplus
\Pic(X)\oplus CH_0(X)\ ,$$where the group $CH_0(X)$ 
parametrizes $0$\tx cycles on $X$. While the structure of the Picard
group $\Pic(X)$ is well understood,  this is not the case
for $CH_0(X)$: if $X$ admits a nonzero holomorphic 2-form, 
 it is an infinite-dimensional
vector space over ${\bf Q}$ ([M], [R]). 
\ind The simplest such surfaces are probably the K3 surfaces, which carry
a nowhere vanishing holomorphic 2-form. In this case $\Pic(X)$ is a
lattice, while $CH_0(X)$ is very large; the following result is therefore 
somewhat  surprising: \smallskip {\bf Theorem}$.-$ {\it The image of the
intersection product}$$\Pic(X)\otimes\Pic(X)\rightarrow CH_0(X)$${\it is
infinite cyclic}. 
\smallskip
\ind In fact we construct a canonical class $\xi^{} _X\in CH_0(X)$
such that any product of divisors is a multiple of $\xi^{} _X$ (Proposition
1). There is another canonical class  in $CH_0(X)$, namely the second
Chern class
$c_2(X)$; it seems plausible that it is always a multiple of $\xi^{} _X$, but
we are able to check this only in some particular cases (Proposition 2).
 \vskip1cm {\bf Proofs}
\smallskip
\ind We work over the complex numbers. A {\it rational curve} on a
surface is  irreducible, but possibly singular. 
\ind The Theorem follows from a slightly more precise statement:
\smallskip
{\bf Proposition 1$.-$} {\it Let $X$ be a  projective $K3$ surface. Then all
points of $X$ which lie on some rational curve have the same class
$\xi^{}_X$ in
$CH_0(X)\ ;$ if $D$ and $E$ are divisors on $X$, we have $$D\cdot
E=n\,\xi^{}_X\quad{\rm with}\ n=\deg(D\cdot E)\ .$$  \smallskip  Proof} :
Let $R$ be a rational curve  on $X$; it 
 is the image of a generically injective map $j:{\bf
P}^1\rightarrow X$. Put $\xi^{}_R=j_*(p)$, where $p$ is an
arbitrary point of ${\bf P}^1$. For any divisor $D$
on $X$, we have in
$CH_0(X)$
$$R\cdot D=j_*j^*D=j_*(n\,p)= n\,\xi^{}_R\ ,\quad \hbox{with}\
n=\deg(R\cdot D)\ .$$ 
\ind Let $S$ be another
rational curve. If  $\deg(R\cdot S)\not=0$, the above
equality applied to $R\cdot S$ gives $\xi^{}_S=\xi^{}_R$ (recall that
$CH_0(X)$ is torsion free [R]). If  $\deg(R\cdot S)=0$,  choose an
ample divisor $H$; by a theorem of Bogomolov and Mumford
[M-M], $H$ is linearly equivalent to a sum of rational
curves. Since $H$ is connected, we can find a chain $R_0,\ldots,R_k$ of
distinct rational curves such that $R_0=R$, $R_k=S$ and $R_{i}\cap
R_{i+1}\not=\emptyset$ for $i=0,\ldots,k-1$. We conclude from the
preceding case that $\xi^{} _R=\xi^{} _{R_1}=\ldots=\xi^{}_S$. \ind Thus
the class 
$\xi^{}_R$ does not depend on the choice of $R$; let us denote it by
$\xi^{}_X$.  We have $R\cdot D =\deg(R\cdot D)\ \xi^{}_X$  for any divisor
$D$  and any rational curve $R$ on $X$. Since the group $\Pic(X)$ is
spanned by the  classes of  rational curves (again by the
Bogomolov-Mumford theorem), the Proposition follows.\cqfd {\bf
Remarks}$.-$ 1) The  result (and the proof) hold  more generally for  any 
surface $X$ such that:
\indp a)  The Picard group of $X$ is
spanned by the  classes of  rational curves;
\indp b) There exists an ample divisor on $X$ which is a sum of rational
curves. 
\ind This is the case when $X$ admits a non-trivial elliptic fibration over
${\bf P}^1$ with a section, or for some particular surfaces like  Fermat
surfaces in ${\bf P}^3$ with degree  prime to $6$ [S]. \ind
2) Let
$A$ be an abelian surface. According to [Bl], the image of the product map
$\Pic(A)\otimes \Pic(A)\rightarrow CH_0(A)$  has finite index, so the
situation looks rather different from  the K3 case. There is however an
 analogue to the Proposition. Let us work for simplicity in
the ${\bf Q}$\tx vector space $CH^{}_{\bf Q}(A):=CH(A)\otimes{\bf Q}$.
Let $\Pic^+(A)$ be the subspace of $\Pic^{}_{\bf Q}(A)$ fixed by the
action of the involution $a\mapsto -a$. We have a direct sum decomposition
$$\Pic^{}_{\bf Q}(A)=\Pic^+(A)\oplus \Pic^{\rm o}_{\bf
Q}(A)\ ,$$so that $\Pic^+(A)$ is canonically isomorphic to
$NS^{}_{\bf Q}(A)$. Now we claim that
{\it the image of the map $\mu :\Pic^+(A)\otimes
\Pic^+(A)\rightarrow CH^{}_{\bf Q}(A)$ is ${\bf Q}\,[0]$}. This is a direct
consequence of the decomposition of 
$CH^{}_{\bf Q}(A)$ described in [B]: let $k$ be an integer $\ge 2$, and
let ${\bf k}$  be the multiplication by $k$ in $A$. 
 We have ${\bf k}^*D=k^2D$ for any element $D$ of  $\Pic^+(A)$, thus
${\bf k}^*\xi^{} =k^4\xi^{} $ for any element $\xi^{} $ in the image of 
$\mu $; but the latter property characterizes the multiples of $[0]$.\cqfd 
\smallskip 
 \ind Proposition 1 provides a distinguished class $\xi^{}_X$ of degree 1 in 
the Chow group $CH_0(X)$: by definition, any $0$\tx cycle of degree $d$
whose support is contained in a finite union of rational curves has class
$d\,\xi^{}_X$ in $CH_0(X)$. On the other hand, there is another canonical
element in that group, namely the Chern class $c_2(X)$; we are led to ask
whether this class is proportional to $\xi^{}_X$. In all cases we were able to
check the answer is positive: \smallskip {\bf Proposition 2}$.-$ {\it The
relation $c_2(X)=24\,\xi^{} _X$  in $CH_0(X)$ holds  in the following
cases:}
\ind a) {\it $X$ is a complete intersection in a product of projective
spaces} ; \ind b) {\it $X$ admits an elliptic fibration} ; \ind c) {\it $X$ is a
Kummer surface}. \smallskip {\it Proof} : a) We observe more generally that
if a projective manifold $V$ is such that $c_2(V)$ is a rational combination
of products of divisors, the same property holds for a smooth hypersurface
$Y$ in $V$: this follows at once from the exact sequence $$0\rightarrow
{\cal O}_Y(-Y)\longrightarrow \Omega^1_{V|Y}\longrightarrow
\Omega^1_Y\rightarrow 0\ .$$ \ind b) Let $f:X\rightarrow {\bf P}^1$ be
an elliptic fibration. We have an exact
sequence
$$0\rightarrow f^*\omega_{{\bf P}^1}(V)\longrightarrow
\Omega^1_X\longrightarrow {\cal
I}_Z\otimes f^*\omega_{{\bf P}^1}^{-1}(-V)
\rightarrow 0\ ,$$where $V$ is a sum of smooth rational
curves contained in the fibres of $f$, and
${\cal I}_Z$ is the ideal sheaf of a finite subscheme $Z\i
X$ contained in the locus where $f$ is not smooth. Let $[Z]$
be the class of the corresponding $0$\tx
cycle in $CH_0(X)$, and $F\in \Pic(X)$ the class of a fibre;
from the exact sequence we obtain 
$\ c_2(X)=[Z]-(2F-V)^2$ . But $Z$ is supported by (smooth)
rational curves, hence is proportional to $\xi^{}_X$.
\ind c) Let $X$ be a Kummer surface. There exist  an abelian
surface $A$ and a diagram\vskip-20pt
$$\diagram{\widehat{A}&\hfl{\varepsilon }{}&A\cr
\vfl{\pi }{}&&\cr X &&}$$\vskip-12pt
where
 $\varepsilon:\widehat{A}\rightarrow A$ is the
blowing up of the points of order $\le 2$ of $A$ and $\pi
:\widehat{A}\rightarrow X$ is the quotient map by the involution of
$\widehat{A}$ deduced from $a\mapsto -a$. Let $E$ be  the exceptional
divisor in $\widehat{A}$.  Using the exact
sequences $$\nospacedmath\displaylines{0\rightarrow \pi
^*\Omega^1_X\longrightarrow \Omega^1_{\widehat{A}}\longrightarrow
{\cal O}_E(-E)\rightarrow 0\cr
0\rightarrow \varepsilon^*\Omega^1_A\longrightarrow
\Omega^1_{\widehat{A}}\longrightarrow \Omega^1_E\rightarrow
0}$$a straightforward computation gives $\pi
^*c_2(X)=-3E^2$. Let $C$ be the divisor $\pi _*E$ (sum of
the 16 $(-2)$\tx curves of the Kummer surface). We have $\pi
^*C=2E$ and therefore $\pi ^*C^2=4E^2$, from which we get
$c_2(X)=-{3\over 4}C^2=24\,\xi^{}_X$.\cqfd

\ind It is tempting to conjecture that the relation
$c_2(X)=24\,\xi^{}_X$ always holds, but the evidence is not
overwhelming: the crucial case would be when $\Pic(X)$ has
rank 1 (in which case the theorem is trivial). 
\ind  There seems to be no
reason to expect that the theorem still holds  for other surfaces, say for
regular surfaces of general type. But finding a counter-example is probably
hard: we know of no way to prove that an explicitely given $0$\tx cycle of
degree 0 on such a surface has nonzero class in the Chow group. 
\vskip2cm
\centerline{ REFERENCES} \vglue15pt\baselineskip12.8pt
\def\num#1{\smallskip\item{\hbox to\parindent{\enskip [#1]\hfill}}}
\parindent=1.3cm 
\num{B} A. {\pc BEAUVILLE}: {\sl 	Sur l'anneau de Chow d'une vari\'et\'e
ab\'elienne.} Math. Annalen {\bf 273} (1986),  647--651.
\num{Bl} S. {\pc BLOCH}: {\sl Some elementary theorems about algebraic
cycles on Abelian varieties}. Invent. Math. {\bf 37} (1976),  215--228. 
\num{M} D. {\pc MUMFORD}: {\sl Rational equivalence of $0$-cycles on
surfaces}.  J. Math. Kyoto Univ. {\bf 9} (1968), 195--204. 
\num{M-M} S. {\pc MORI}, S. {\pc MUKAI}: {\sl Mumford's theorem on
curves on $K3$ surfaces}. Algebraic Geometry
(Tokyo/Kyoto 1982), LN {\bf 1016},  351--352; Springer-Verlag (1983).
\num{R} A. A. {\pc ROJTMAN}: {\sl The torsion of the group of $0$-cycles 
modulo rational equi\-valence}.  Ann. of Math. {\bf 111} (1980),
553--569. 
\num{S} T. {\pc SHIODA}: {\sl On the Picard number of a Fermat surface}. J.
Fac. Sci. Univ. Tokyo {\bf 28} (1982), 725--734.
\vskip1cm
\def\pc#1{\eightrm#1\sixrm}
\hfill\vtop{\eightrm\hbox to 5cm{\hfill Arnaud {\pc BEAUVILLE}\hfill}
 \hbox to 5cm{\hfill Laboratoire J.-A. Dieudonn\'e\hfill}
 \hbox to 5cm{\sixrm\hfill UMR 6621 du CNRS\hfill}
\hbox to 5cm{\hfill {\pc UNIVERSIT\'E DE}  {\pc NICE}\hfill}
\hbox to 5cm{\hfill  Parc Valrose\hfill}
\hbox to 5cm{\hfill F-06108 {\pc NICE} Cedex 02\hfill}}
\end